\documentclass[12pt,reqno]{amsart}

\usepackage[utf8]{inputenc}
\usepackage[T1]{fontenc}
\usepackage{lmodern}
\usepackage{amsfonts,amssymb,mathrsfs}
\usepackage{microtype}
\usepackage[margin=1in]{geometry}
\usepackage[colorlinks=true,allcolors=blue]{hyperref}
\makeatletter
\long\def\@makefntext#1{\noindent #1}
\makeatother

\numberwithin{equation}{section}
\theoremstyle{plain}
\newtheorem{theorem}{Theorem}[section]
\newtheorem{lemma}{Lemma}[section]
\newtheorem{proposition}{Proposition}[section]
\newtheorem{corollary}{Corollary}[section]
\theoremstyle{definition}
\newtheorem{definition}{Definition}[section]
\theoremstyle{remark}
\newtheorem{remark}{Remark}[section]

\begin{document}

\title[Lower bounds for Kolmogorov widths]{Lower bounds for the Kolmogorov widths of global attractors of dissipative systems}
\author{Arthur Sinai}
\thanks{Yaroslav-the-Wise Novgorod State University, Veliky Novgorod, Russia.\\
E-mail: \texttt{arthur.sinai@mail.ru}.}
\date{}

\begin{abstract}
We introduce the notion of a strengthening set in a Hilbert space and prove an exponential lower bound for its Kolmogorov widths. For compact sets containing the origin, the strengthening property is equivalent to the existence of such a bound. As an application, we obtain an exact formula for the widths of the global attractor of a diagonal dissipative system.
\end{abstract}

\maketitle

\section{Introduction}\label{sec-intro}

The notion of Kolmogorov width was introduced by A.\,N.~Kolmogorov~\cite{Kol36}. Width theory in function spaces is surveyed by V.\,M.~Tikhomirov~\cite{Tik60}, and a general treatment is given in~\cite{Pin85}. The Kolmogorov width of order \(n\) of a set \(M\) in a normed space \(X\) is the quantity
\[
d_n(M,X)=\inf_{L\subset X}\sup_{x\in M}\inf_{y\in L}\|x-y\|_X,
\]
where the infimum is taken over all linear subspaces \(L\) of \(X\) with dimension at most \(n\). The width gives the best uniform error of approximation of \(M\) by such subspaces. For a totally bounded set \(M\), its \(\varepsilon\)-entropy with respect to \(X\) is the base-two logarithm of the smallest number of closed balls of radius \(\varepsilon\) in \(X\) covering \(M\)~\cite{KolTik59}.

Bounds for the fractal and Hausdorff dimensions of global attractors of dissipative systems have been studied extensively~\cite{BabVis83,Tem97,CheVis02,Zel23}. Eden, Zelik and Kalantarov~\cite[\S~5]{EZK13} constructed an attractor of a semilinear parabolic system whose orthogonal projection contains line segments parallel to the vectors of an orthonormal system. Below we derive a lower bound for its widths from this construction.

An upper bound for the widths of an attractor can often be obtained by choosing a subspace adapted to the smoothness of the attractor's elements. For example, consider the orthonormal Fourier system on the torus \(\mathbb T^d\). If, for some \(\rho>0\),
\[
B=\sup_{u\in\mathcal A}\sum_{k\in\mathbb Z^d}
e^{2\rho|k|}|\widehat u(k)|^2<\infty,
\]
then truncating the Fourier series to the indices \(|k|\leqslant R\) gives an error in \(L^2(\mathbb T^d)\) of at most \(\sqrt B\,e^{-\rho R}\), using a subspace of dimension of order \(R^d\). Hence \(d_n(\mathcal A)\leqslant C_{\mathrm{up}}e^{-c n^{1/d}}\). To obtain a lower bound, the approximation error must be bounded from below for every subspace of dimension at most \(n\).

A.~Cohen~\cite{Coh14} explicitly identified the problem of obtaining lower bounds for the widths of solution sets of parametric partial differential equations as open.

Our approach is based on the notion of a \emph{strengthening set} (Definition~\ref{def-spread}). Let \(X\) be a Hilbert space. A bounded set \(A\subset X\) with \(0\in A\) is called a strengthening set with parameters \((C,\kappa)\), \(C,\kappa>0\), if there exists an orthonormal system \(\{e_j\}_{j\geqslant0}\) such that, for every \(j\geqslant0\) and every finite-dimensional subspace \(L\subset X\), there exists \(x\in A\) such that
\[
\inf_{y\in L}\|x-y\|
\geqslant Ce^{-\kappa j}\inf_{z\in L}\|e_j-z\|.
\]
A weighted bound based on Bessel's inequality (Theorem~\ref{th-l2}) gives the lower bound \(d_n(A)\geqslant C a(\kappa)e^{-\kappa n}\) for a strengthening set, with the same exponent. The constant \(a(\kappa)\) is determined exactly. A discrete ellipsoid has the smallest widths at each \(n\) among sets with prescribed strengthening parameters (Theorem~\ref{th-rootfree} and Proposition~\ref{prop-ellipsoid_sharp}).

For a compact set \(A\) with \(0\in A\), the strengthening property with prescribed parameters holds if and only if its closed absolutely convex hull contains line segments \([-r_je_j,r_je_j]\) along an orthonormal system, with half-lengths \(r_j\) satisfying an exponential lower bound with the same parameters (Lemma~\ref{lem-absconv_inner}). To deduce the strengthening property from a lower bound for the widths, we use a bound for sections of this hull. From the inequality \(d_n(A)\geqslant Ce^{-\kappa n}\) for all \(n\) we construct an orthonormal system for which \(r_j\) can be chosen so that \(r_j\geqslant Ce^{-\kappa j}/(1+\sqrt j)\). In particular, \(A\) is a strengthening set with exponent \(\kappa+\varepsilon\) for every \(\varepsilon>0\) (Theorem~\ref{th-main_dichotomy}). Thus the strengthening property is equivalent to the existence of an exponential lower bound for the widths (Corollary~\ref{cor-dichotomy_equiv}).

Orthogonal spread gives a simple sufficient condition for the strengthening property (Proposition~\ref{prop-orth_to_spread}). Under almost orthogonal spread, square summability of the relative remainders allows the exponent in the lower bound to be preserved (Theorem~\ref{th-almost_summable}).

For a diagonal dissipative system we explicitly determine the global attractor and obtain an exact formula for its widths. A lower bound for the widths in terms of the expansion coefficients of stationary solutions with square-summable relative remainders is given in Proposition~\ref{prop-stationary}.

\textbf{Organization of the paper.}
Section~\ref{sec-prelim} contains the preliminaries.
The main lower bounds for the widths are proved in Section~\ref{sec-main},
and sufficient conditions and examples are given in Section~\ref{sec-examples}.
Section~\ref{sec-dichotomy} is devoted to the converse theorem and the dichotomy,
and Section~\ref{sec-appl} to applications to dissipative systems.

\section{Preliminaries}\label{sec-prelim}

Throughout, \(X\) is a real Hilbert space with inner product \(\langle \cdot, \cdot \rangle\) and norm \(\|\cdot\|\). Put \(\mathbb Z_+=\{0,1,\ldots\}\) and \(\mathbb N=\{1,2,\ldots\}\).

\begin{definition}\label{def-kolmogorov}
\begin{sloppypar}
The Kolmogorov width of order \(n\) of a set \(M\subset X\) is the quantity
\[
d_n(M,X)=\inf_{L\subset X}\sup_{x\in M}\inf_{y\in L}\|x-y\|_X,
\]
where the infimum is taken over all linear subspaces \(L\) of \(X\) with dimension at most \(n\).
\end{sloppypar}
\end{definition}

When \(X\) is fixed, we write \(d_n(M)\) instead of \(d_n(M,X)\).

\begin{lemma}\label{lem-infDist_smul}
For a linear subspace \(L \subset X\), a vector \(v \in X\) and a scalar \(c \in \mathbb{R}\), we have
\[
\inf_{y \in L} \|cv - y\| = |c| \inf_{z \in L} \|v - z\|.
\]
\end{lemma}

\begin{proof}
For \(c=0\) the equality is immediate; for \(c\ne0\) it follows from the substitution \(y=cz\) and the linearity of \(L\).
\end{proof}

\begin{lemma}\label{lem-pigeonhole}
Let \(\{e_0, \dots, e_m\}\) be an orthonormal system in \(X\), and let \(L \subset X\) be a linear subspace with \(\dim L \leqslant m\). Then there is an index \(j \in \{0, \dots, m\}\) such that
\[
\inf_{y \in L} \|e_j - y\| \geqslant \frac{1}{\sqrt{m+1}}.
\]
\end{lemma}

\begin{proof}
Let \(P_L\) be the orthogonal projection onto \(L\), and let \(\{b_i\}_{i=1}^{s}\) be an orthonormal basis of \(L\). By Bessel's inequality,
\[
\sum_{j=0}^{m}\|P_L e_j\|^2
=\sum_{i=1}^{s}\sum_{j=0}^{m}|\langle e_j,b_i\rangle|^2
\leqslant s=\dim L.
\]
Thus the sum of the squared distances from \(e_j\) to \(L\) is at least \(m+1-\dim L\geqslant1\). At least one of these distances is not less than \(1/\sqrt{m+1}\).
\end{proof}

\begin{lemma}\label{lem-weighted_pigeonhole}
Let \(\{e_0, \dots, e_{M-1}\}\) be an orthonormal system in \(X\), let \(L \subset X\) be a linear subspace with \(\dim L \leqslant n < M\), and let \(r_0, \dots, r_{M-1} > 0\) be weights. Then there is an index \(j \in \{0, \dots, M-1\}\) such that
\[
r_j \inf_{z \in L} \|e_j - z\| \geqslant \Bigl(\frac{M-n}{\sum_{i=0}^{M-1} r_i^{-2}}\Bigr)^{1/2}.
\]
For \(M=n+1\) and equal weights \(r_i \equiv 1\), the right-hand side is \(1/\sqrt{n+1}\), and the lemma reduces to Lemma~\ref{lem-pigeonhole}.
\end{lemma}

\begin{proof}
Put \(h_j=\inf_{z\in L}\|e_j-z\|\) and \(D=\max_{j<M} r_jh_j\). As in the proof of Lemma~\ref{lem-pigeonhole},
\[
M-n\leqslant\sum_{j<M}h_j^2\leqslant D^2\sum_{j<M}r_j^{-2}.
\]
This gives the required bound for an index at which the maximum is attained.
\end{proof}

\section{Strengthening sets and the main lower bound}\label{sec-main}

\begin{definition}\label{def-spread}
A bounded set \(A \subset X\) with \(0 \in A\) is called a \emph{strengthening set} with parameters \((C, \kappa)\), \(C, \kappa > 0\), if there exist vectors \(e_k \in X\) and numbers \(r_k\) such that:
\begin{enumerate}
\item[(i)] \(\{e_k\}_{k=0}^{\infty}\) is an orthonormal system;
\item[(ii)] \(r_k \geqslant C e^{-\kappa k}\) for all \(k \geqslant 0\);
\item[(iii)] for every \(j \geqslant 0\) and every finite-dimensional subspace \(L \subset X\) there is a point \(x \in A\) such that
\[
\inf_{y \in L} \|x - y\| \geqslant r_j \inf_{z \in L} \|e_j - z\|.
\]
\end{enumerate}
\end{definition}

To bound \(d_n(A)\), it is enough to assume condition~(iii) for finitely many directions.

\begin{theorem}\label{th-l2}
Let \(A \subset X\), let \(\{e_0, \dots, e_{M-1}\}\) be an orthonormal system and let \(r_0, \dots, r_{M-1} > 0\). Suppose that condition~(iii) of Definition~\ref{def-spread} holds for \(j = 0, \dots, M-1\). Then, for \(0 \leqslant n < M\),
\begin{equation}\label{eq-l2M}
d_n(A) \geqslant \Bigl(\frac{M-n}{\sum_{j=0}^{M-1} r_j^{-2}}\Bigr)^{1/2}.
\end{equation}
In particular, if \(A\) is a strengthening set with an orthonormal system \(\{e_k\}\) and numbers \(r_k\), then, for every \(n \in \mathbb{Z}_+\),
\begin{equation}\label{eq-l2}
d_n(A) \geqslant \mu_n, \qquad \mu_n := \Bigl(\sum_{j=0}^{n} r_j^{-2}\Bigr)^{-1/2}.
\end{equation}
\end{theorem}

\begin{proof}
For any \(L\) of dimension at most \(n\), condition~(iii) gives
\[
\sup_{x\in A}\inf_{y\in L}\|x-y\|
\geqslant \max_{j<M}r_j\inf_{z\in L}\|e_j-z\|.
\]
Applying Lemma~\ref{lem-weighted_pigeonhole} to this maximum and taking the infimum over \(L\), we obtain~\eqref{eq-l2M}. The bound~\eqref{eq-l2} is the case \(M=n+1\).
\end{proof}

For the set \(\{r_je_j:0\leqslant j<M\}\), the bound~\eqref{eq-l2M} also follows from Ismagilov's theorem~\cite[Theorem~1]{Ism68} by taking a probability measure with masses proportional to \(r_j^{-2}\).

\begin{remark}\label{rem-fixed_average}
The strengthening method can give a lower bound with a constant independent of the dimension, while averaging with fixed weights may lose this uniformity. Consider
\[
B=\{0\}\cup\{r_je_j:j\geqslant0\},\qquad r_j=Ce^{-\kappa j},
\]
where \(\{e_j\}\) is an orthonormal system and \(C,\kappa>0\). By Theorem~\ref{th-l2},
\[
d_n(B)\geqslant\sqrt{1-e^{-2\kappa}}\,r_n.
\]
When the squared errors are averaged with arbitrary fixed weights \(p_j\geqslant0\), \(\sum_jp_j=1\), choosing the subspace spanned by the first \(n\) coordinate vectors gives
\[
\sum_{j\geqslant n}p_jr_j^2
\leqslant r_n^2\sum_{j\geqslant n}p_j
=o(r_n^2).
\]
Thus such averaging does not yield a lower bound of order \(r_n\) with a positive constant common to all \(n\). If the weights are allowed to depend on \(n\), averaging reproduces the main bound~\eqref{eq-l2M}.
\end{remark}

\begin{theorem}\label{th-rootfree}
Let \(A\) be a strengthening set with parameters \((C,\kappa)\), and let \(q=e^{-\kappa}\). Then, for all \(n\in\mathbb Z_+\),
\begin{equation}\label{eq-rootfree}
d_n(A)\geqslant C q^n
\left((1-q^2)\max_{p\in\mathbb N}
\frac{p q^{2(p-1)}}{1-q^{2(n+p)}}\right)^{1/2}
\geqslant C a(\kappa) q^n,
\end{equation}
where
\[
a(\kappa)=\left((1-q^2)\max_{p\in\mathbb N}p q^{2(p-1)}\right)^{1/2}.
\]
Within the class of strengthening sets with parameters \((C,\kappa)\), the first bound is attained at each \(n\), and the constant \(a(\kappa)\) in the second bound is optimal uniformly in \(n\). The maximum in the formula for \(a(\kappa)\) is attained at \(p=\lfloor(1-q^2)^{-1}\rfloor\).
\end{theorem}

\begin{proof}
Apply Theorem~\ref{th-l2} with \(M=n+p\). Since \(r_j\geqslant Cq^j\), we have
\[
d_n(A)^2\geqslant\frac{p}{\sum_{j<n+p}r_j^{-2}}
\geqslant C^2 q^{2n}\frac{(1-q^2)pq^{2(p-1)}}{1-q^{2(n+p)}}.
\]
Taking the maximum over \(p\) gives the first bound, and replacing the denominator by one gives the second. The ratio of consecutive terms \(p q^{2(p-1)}\) is \(q^2(p+1)/p\), which gives the stated choice of \(p\). Optimality follows from Proposition~\ref{prop-ellipsoid_sharp}.
\end{proof}

For \(\kappa\geqslant\tfrac12\ln2\) we have \(a(\kappa)=\sqrt{1-e^{-2\kappa}}\). As \(\kappa\downarrow0\), taking \(p=\lfloor(1-q^2)^{-1}\rfloor\) gives \((1-q^2)p\to1\) and \(q^{2(p-1)}\to e^{-1}\), whence \(a(\kappa)\to e^{-1/2}\).

\subsection{\texorpdfstring{Sharpness of the \(\ell^2\)-bound}{Sharpness of the l2 bound}}\label{sec-l2_sharp}

\begin{theorem}\label{th-l2_sharp}
Let \(\{e_j\}_{j=0}^{n}\) be an orthonormal system, let \(r_0,\dots,r_n>0\), and put \(F_n = \{0\}\cup\{r_j e_j : 0\leqslant j\leqslant n\}\). Then
\[
d_n(F_n) = \mu_n = \Bigl(\sum_{j=0}^{n} r_j^{-2}\Bigr)^{-1/2}.
\]
\end{theorem}

\begin{proof}
The lower bound follows from Theorem~\ref{th-l2} with \(M=n+1\), since \(F_n\) contains all the points \(r_j e_j\), \(j\leqslant n\). In the space \(\operatorname{span}\{e_0,\ldots,e_n\}\) take the hyperplane with unit normal
\[
w=\sum_{j=0}^{n}\frac{\mu_n}{r_j}e_j.
\]
The distance from \(r_j e_j\) to this hyperplane is \(|\langle r_j e_j,w\rangle|=\mu_n\) for all \(j\leqslant n\).
\end{proof}

\subsection{Other rates of decay}\label{sec-decay}

Theorem~\ref{th-l2} does not use condition~(ii) of Definition~\ref{def-spread} and applies to arbitrary positive \(r_k\).

\begin{corollary}\label{cor-decay_poly}
Suppose that, for a set \(A\subset X\), condition~(iii) of Definition~\ref{def-spread} holds for all \(j\geqslant0\) with an orthonormal system \(\{e_k\}_{k\geqslant0}\) and \(r_k=(k+1)^{-\alpha}\), \(\alpha\geqslant0\). Then
\[
d_n(A)\geqslant\frac{(2n+1)^{-\alpha}}{\sqrt2}\geqslant2^{-\alpha-1/2}(n+1)^{-\alpha}, \qquad n\in\mathbb{Z}_+.
\]
The order \((n+1)^{-\alpha}\) is sharp: for \(A=\{0\}\cup\{r_ke_k:k\geqslant0\}\) we have \(d_n(A)\leqslant(n+1)^{-\alpha}\).
\end{corollary}

\begin{proof}
Apply Theorem~\ref{th-l2} with \(M=2n+1\). Since \(r_j\geqslant(2n+1)^{-\alpha}\) for \(j\leqslant2n\),
\[
\sum_{j=0}^{2n} r_j^{-2} \leqslant (2n+1)^{2\alpha+1}, \qquad d_n(A)\geqslant\Bigl(\frac{n+1}{(2n+1)^{2\alpha+1}}\Bigr)^{1/2}\geqslant\frac{(2n+1)^{-\alpha}}{\sqrt2}.
\]
The second inequality in the statement follows from \(2n+1\leqslant2(n+1)\). For the specified set \(A\), it suffices to take \(x=r_je_j\) in condition~(iii). The subspace \(\operatorname{span}\{e_0,\dots,e_{n-1}\}\) approximates every point of \(A\) with error at most \(r_n=(n+1)^{-\alpha}\).
\end{proof}

For the set \(\{0\}\cup\{(k+1)^{-\alpha}e_k\}\) we have \(d_n\asymp(n+1)^{-\alpha}\), whereas \(\mu_n\asymp(n+1)^{-\alpha-1/2}\). Passing from \(n+1\) to \(2n+1\) directions gives a lower bound of the correct order.

\begin{corollary}\label{cor-decay_lacunary}
Suppose that, for a set \(A\subset X\), condition~(iii) of Definition~\ref{def-spread} holds for all \(j\geqslant0\) with an orthonormal system \(\{e_k\}_{k\geqslant0}\) and positive \(r_k\). If \(\sum_{j=0}^n r_j^{-2}\leqslant K\,r_n^{-2}\) for all \(n\) with a single constant \(K>0\), then \(d_n(A)\geqslant r_n/\sqrt K\).
\end{corollary}

\begin{proof}
By Theorem~\ref{th-l2},
\[
d_n(A)\geqslant\mu_n = \Bigl(\sum_{j=0}^n r_j^{-2}\Bigr)^{-1/2} \geqslant \frac{r_n}{\sqrt K}.
\]
\end{proof}

For the geometric sequence \(r_k=Ce^{-\kappa k}\), the hypothesis of the corollary holds with \(K=(1-e^{-2\kappa})^{-1}\) by the formula for the sum of a geometric progression. This gives the bound~\eqref{eq-rootfree} with \(p=1\) after replacing the denominator by one.

The hypothesis of Corollary~\ref{cor-decay_lacunary} is sufficient but not necessary for the bound \(d_n(A)\geqslant c r_n\). For \(r_k=(k+1)^{-\alpha}\), Corollary~\ref{cor-decay_poly} gives \(d_n(A)\geqslant 2^{-\alpha-1/2}r_n\), although this hypothesis fails.

\section{Sufficient conditions and examples}\label{sec-examples}

\begin{definition}\label{def-orth}
A set \(A \subset X\) with \(0 \in A\) is said to have \emph{orthogonal spread} with parameters \((C, \kappa)\), \(C,\kappa>0\), if there exist an orthonormal system \(\{e_k\}_{k=0}^{\infty}\) and points \(x_k = c_k e_k \in A\) with \(|c_k| \geqslant C e^{-\kappa k}\).
\end{definition}

\begin{proposition}\label{prop-orth_to_spread}
Every bounded set with orthogonal spread is a strengthening set with the same parameters.
\end{proposition}

\begin{proof}
Put \(r_k = |c_k|\). For any \(j \geqslant 0\) and \(L \subset X\), choose \(x = x_j \in A\). By Lemma~\ref{lem-infDist_smul},
\[
\inf_{y \in L} \|x_j - y\| = |c_j| \inf_{z \in L} \|e_j - z\| = r_j \inf_{z \in L} \|e_j - z\|.
\]
Condition~(iii) of Definition~\ref{def-spread} holds.
\end{proof}

\begin{corollary}\label{cor-orth_bound}
If \(A\) is bounded and has orthogonal spread with parameters \((C, \kappa)\), then the following lower bound holds with the same exponential rate:
\[
d_n(A) \geqslant C a(\kappa)e^{-\kappa n}.
\]
\end{corollary}

\begin{proof}
By Proposition~\ref{prop-orth_to_spread}, the set \(A\) is a strengthening set with the same parameters, so the bound follows from Theorem~\ref{th-rootfree}.
\end{proof}

\begin{definition}\label{def-almost}
A set \(A \subset X\) with \(0\in A\) is said to have \emph{almost orthogonal spread} with parameters \((C,\kappa,\delta,(\eta_k))\), where \(C,\kappa>0\), \(\delta\in(0,1]\) and \(\eta_k\geqslant0\), if there exist an orthonormal system \(\{e_k\}_{k\geqslant0}\) and points \(x_k\in A\) for which the following conditions hold for all \(k\geqslant0\):
\begin{enumerate}
\item[(i)] \(|\langle x_k, e_k \rangle| \geqslant \delta \|x_k\|\);
\item[(ii)] \(\|x_k\| \geqslant C e^{-\kappa k}\);
\item[(iii)] \(\bigl\|x_k - \langle x_k, e_k \rangle e_k\bigr\|^2 \leqslant \eta_k |\langle x_k, e_k \rangle|^2\).
\end{enumerate}
\end{definition}

\begin{theorem}\label{th-almost_bound}
Suppose that \(A\) is bounded and has almost orthogonal spread with parameters \((C,\kappa,\delta,(\eta_k))\). If \(m,n\in\mathbb Z_+\) and \(4 \hat{\eta}_{m,n} (n+1) < 1\), where \(\hat{\eta}_{m,n} = \max_{k \in \{m, \dots, m+n\}} \eta_k\), then
\[
d_n(A) \geqslant \frac{\delta}{2} \frac{C e^{-\kappa (m+n)}}{\sqrt{n+1}}.
\]
\end{theorem}

\begin{proof}
Let \(\dim L\leqslant n\). By Lemma~\ref{lem-pigeonhole}, there exists \(j\in\{m,\ldots,m+n\}\) such that \(\|P_{L^\perp}e_j\|\geqslant1/\sqrt{n+1}\). Put \(a_j=\langle x_j,e_j\rangle\) and \(w_j=x_j-a_je_j\). Then \(\|w_j\|\leqslant\sqrt{\hat\eta_{m,n}}|a_j|\), and
\[
\|P_{L^\perp}x_j\|
\geqslant |a_j|\left(\frac{1}{\sqrt{n+1}}-\sqrt{\hat\eta_{m,n}}\right)
\geqslant\frac{\delta C e^{-\kappa(m+n)}}{2\sqrt{n+1}}.
\]
In the last bound we used the condition \(4\hat\eta_{m,n}(n+1)<1\) and the inequality \(|a_j|\geqslant\delta C e^{-\kappa j}\). Taking the infimum over \(L\) proves the assertion.
\end{proof}

If \(\sum_k\eta_k<\infty\), the exponent \(\kappa\) can be preserved and the factor \(1/\sqrt{n+1}\) can be removed from the lower bound.

\begin{theorem}\label{th-almost_summable}
Suppose that \(A\) is bounded and has almost orthogonal spread with parameters \((C,\kappa,\delta,(\eta_k))\), and that \(\sum_{k\geqslant0}\eta_k<\infty\). If \(m \in \mathbb{Z}_+\) is such that \(\sum_{k\geqslant m}\eta_k\leqslant\frac14\), then, for all \(n \in \mathbb{Z}_+\),
\begin{equation}\label{eq-almost_summable}
d_n(A) \geqslant \frac{\delta}{2}\, C e^{-\kappa m}a(\kappa)e^{-\kappa n}.
\end{equation}
\end{theorem}

\begin{proof}
Put \(a_k=\langle x_k,e_k\rangle\) and \(w_k=x_k-a_ke_k\). By conditions~(i)--(iii) of Definition~\ref{def-almost}, we have \(|a_k|\geqslant\delta Ce^{-\kappa k}>0\) and \(\|w_k\|/|a_k|\leqslant\sqrt{\eta_k}\). Define
\[
Tx=x+\sum_{k\geqslant m}\langle x,e_k\rangle\,\frac{w_k}{a_k}.
\]
By the Cauchy--Schwarz and Bessel inequalities, the series converges absolutely and \(\|T-I\|\leqslant(\sum_{k\geqslant m}\eta_k)^{1/2}\leqslant1/2\). Hence \(T\) is invertible and \(\|T^{-1}\|\leqslant2\). For \(k\geqslant m\) we have \(T(a_ke_k)=x_k\).

The set \(B=\{0\}\cup\{a_ke_k:k\geqslant m\}\) is bounded, since \(|a_k|\leqslant\|x_k\|\). It has orthogonal spread with parameters \((\delta Ce^{-\kappa m},\kappa)\) relative to the orthonormal system \(\{e_{m+j}\}_{j\geqslant0}\), and \(T(B)\subset A\). For any subspace \(L\) of dimension at most \(n\),
\[
\inf_{y\in L}\|Tx-y\|
\geqslant \tfrac12\inf_{z\in T^{-1}L}\|x-z\|.
\]
It follows that \(d_n(A)\geqslant d_n(T(B))\geqslant\tfrac12d_n(B)\). Now apply Corollary~\ref{cor-orth_bound} to \(B\).
\end{proof}

The hypothesis of Theorem~\ref{th-almost_bound} may also hold when \(\sum_k\eta_k=\infty\). For example, for \(\eta_k=(k+1)^{-1}\) it suffices to take \(m=4(n+1)\).

\begin{definition}\label{def-box}
Let \(\{e_k\}_{k\geqslant0}\) be an orthonormal system, and let \(r_k=Ce^{-\kappa k}\), where \(C,\kappa>0\). We call the set
\[
A_{\mathrm{par}} = \{0\} \cup \Bigl\{\sum_{k=0}^{\infty} \varepsilon_k r_k e_k \mathrel{:} \varepsilon_k \in \{-1, 1\}\Bigr\}
\]
a \emph{parallelepiped set} (the series converges in \(X\) because \(\sum r_k^2 < \infty\)).
\end{definition}

\begin{lemma}\label{lem-box_spread}
The set \(A_{\mathrm{par}}\) is a strengthening set with parameters \((C,\kappa)\).
\end{lemma}

\begin{proof}
The set \(A_{\mathrm{par}}\) contains the origin and is bounded, since every nonzero element has norm \(\bigl(\sum_k r_k^2\bigr)^{1/2}\). For any \(j\geqslant0\) and any finite-dimensional \(L\subset X\), take two points \(x_+,x_-\in A_{\mathrm{par}}\) whose signs agree at all indices except \(j\). Then \(x_+-x_-=2r_je_j\), and the triangle inequality gives
\[
\max\{\|P_{L^\perp}x_+\|,\|P_{L^\perp}x_-\|\}
\geqslant r_j\|P_{L^\perp}e_j\|.
\]
Thus condition~(iii) of Definition~\ref{def-spread} holds.
\end{proof}

\begin{proposition}\label{prop-box_exact}
For all \(n \in \mathbb{Z}_+\),
\[
d_n(A_{\mathrm{par}}) = \Bigl(\sum_{k\geqslant n} r_k^2\Bigr)^{1/2} = \frac{C e^{-\kappa n}}{\sqrt{1 - e^{-2\kappa}}}.
\]
In particular, \(d_n(A_{\mathrm{par}}) > r_n \geqslant \mu_n\).
\end{proposition}

\begin{proof}
For \(L_n = \operatorname{span}\{e_0, \dots, e_{n-1}\}\) and any nonzero \(x \in A_{\mathrm{par}}\) we have \(\inf_{y \in L_n} \|x - y\|^2 = \sum_{k\geqslant n} r_k^2\), which gives the upper bound. For the lower bound, fix \(L\) of dimension at most \(n\) and \(N>n\). Put \(p_k=\|P_L e_k\|^2\in[0,1]\). By Bessel's inequality, \(\sum_k p_k \leqslant n\). Average \(\|P_{L^\perp} x\|^2\) over all \(2^N\) points \(x = \sum_k \varepsilon_k r_k e_k\) with \(\varepsilon_k = 1\) for \(k \geqslant N\). The cross terms containing \(\varepsilon_k\) with \(k < N\) vanish under averaging, so the average is at least
\[
\sum_{k<N} r_k^2 \|P_{L^\perp} e_k\|^2 = \sum_{k<N} r_k^2 (1 - p_k).
\]
Since the sequence \(r_k\) is decreasing,
\[
\sum_{k<N} r_k^2 p_k - \sum_{k<n} r_k^2 = \sum_{k<n} r_k^2 (p_k - 1) + \sum_{n\leqslant k<N} r_k^2 p_k \leqslant r_n^2 \Bigl(\sum_{k<N} p_k - n\Bigr) \leqslant 0.
\]
Thus \(\sup_{x \in A_{\mathrm{par}}} \inf_{y \in L} \|x - y\|^2 \geqslant \sum_{n\leqslant k<N} r_k^2\), and letting \(N \to \infty\) gives the lower bound. Finally, \(r_n < \bigl(\sum_{k\geqslant n} r_k^2\bigr)^{1/2}\) and \(\mu_n \leqslant (r_n^{-2})^{-1/2} = r_n\).
\end{proof}

\begin{proposition}\label{prop-box_no_almost}
The set \(A_{\mathrm{par}}\) does not have almost orthogonal spread.
\end{proposition}

\begin{proof}
Let \(\{f_k\}\) be an arbitrary orthonormal system. For every \(k\),
\[
\sup_{x\in A_{\mathrm{par}}}|\langle x,f_k\rangle|
\leqslant\sum_{i\geqslant0} r_i|\langle e_i,f_k\rangle|
\longrightarrow 0\quad(k\to\infty).
\]
Indeed, for each \(i\) we have \(\langle e_i,f_k\rangle\to0\) by Bessel's inequality, and \(\sum_i r_i<\infty\). On the other hand, condition~(i) of Definition~\ref{def-almost} requires \(|\langle x_k,f_k\rangle|\) to be at least \(\delta(\sum_i r_i^2)^{1/2}>0\), since condition~(ii) excludes \(x_k=0\).
\end{proof}

\subsection{A discrete ellipsoid}\label{sec-ellipsoid}
\begin{proposition}\label{prop-ellipsoid_sharp}
Let \(\{e_j\}_{j\geqslant0}\) be an orthonormal system, and let \(\sigma_j>0\) decrease monotonically to zero. For the discrete ellipsoid \(E=\{0\}\cup\{\sigma_j e_j:j\geqslant0\}\),
\begin{equation}\label{eq-ellipsoid_exact}
d_n(E)^2=\max_{m>n}\frac{m-n}{\sum_{j<m}\sigma_j^{-2}}.
\end{equation}
Here the maximum is taken over the integers \(m>n\). In particular, when \(\sigma_j=Cq^j\), the first lower bound in~\eqref{eq-rootfree} is attained, and \(\lim_{n\to\infty}q^{-n}d_n(E)=Ca(\kappa)\), where \(q=e^{-\kappa}\).
\end{proposition}

\begin{proof}
Put \(E_M=\{0\}\cup\{\sigma_j e_j:0\leqslant j<M\}\). The distance to a linear subspace is a seminorm, so \(E_M\) and the octahedron \(\operatorname{absconv}(E_M)\) have the same widths. Projection onto \(\operatorname{span}\{e_0,\ldots,e_{M-1}\}\) does not increase the approximation error, so the widths of \(E_M\) in this span and in \(X\) coincide. For \(n<M\), Sofman's formula for a finite-dimensional octahedron~\cite[formula~(1)]{Sof69} gives the right-hand side of~\eqref{eq-ellipsoid_exact}, with the maximum taken over \(n<m\leqslant M\). On the other hand,
\[
d_n(E_M)\leqslant d_n(E)\leqslant \max\{d_n(E_M),\sigma_M\},
\]
since the distance from any omitted point to any linear subspace is at most \(\sigma_M\). Letting \(M\to\infty\) gives~\eqref{eq-ellipsoid_exact}. The maximum is attained, since for fixed \(n\) the fraction tends to zero as \(m\) tends to infinity.

Substituting \(\sigma_j=Cq^j\) and \(m=n+p\) gives equality in the first bound in~\eqref{eq-rootfree}. Since \(1-q^{2(n+1)}\leqslant 1-q^{2(n+p)}<1\) uniformly in \(p\geqslant1\), the formula for the limit follows as well.
\end{proof}

\section{A dichotomy}\label{sec-dichotomy}

\begin{lemma}\label{lem-section_width}
Let \(K\subset X\) be a compact convex set symmetric about the origin, and let \(L\subset X\) be a subspace of dimension at most \(n\). Then
\[
\max_{z\in K\cap L^\perp}\|z\|\geqslant\frac{d_n(K)}{1+\sqrt n}.
\]
\end{lemma}

\begin{proof}
Put \(\rho=\max_{z\in K\cap L^\perp}\|z\|\) and \(V=\operatorname{span}(P_LK)\). If \(m=\dim V=0\), then \(d_n(K)\leqslant d_0(K)=\rho\). For \(m>0\) the set \(P_LK\) contains the absolutely convex hull of some basis of \(V\). Choose a Euclidean norm \(|\cdot|_*\) on \(V\) in which the L\"owner ellipsoid of \(P_LK\) is the unit ball. By John's theorem~\cite[Theorem~2.1]{Henk12}, there exist contact points \(y_i\in P_LK\) and weights \(a_i>0\) such that
\[
\sum_i a_i\langle y_i,y\rangle_*y_i=y\quad(y\in V),\qquad \sum_i a_i=m.
\]
Choose \(x_i\in K\) with \(P_Lx_i=y_i\) and put \(Ty=\sum_i a_i\langle y_i,y\rangle_*x_i\). Then \(P_LTy=y\), and, for \(y\in P_LK\),
\[
\sum_i a_i|\langle y_i,y\rangle_*|
\leqslant\sqrt m\left(\sum_i a_i|\langle y_i,y\rangle_*|^2\right)^{1/2}
=\sqrt m\,|y|_*\leqslant\sqrt m.
\]
By the absolute convexity of \(K\), for every \(x\in K\) the vector \((x-TP_Lx)/(1+\sqrt m)\) belongs to \(K\cap L^\perp\). Since \(\dim T(V)\leqslant m\leqslant n\), we obtain \(d_n(K)\leqslant(1+\sqrt m)\rho\leqslant(1+\sqrt n)\rho\).
\end{proof}

\begin{remark}\label{rem-gelfand}
For the Gelfand width \(d^n(K)=\inf_{\dim L\leqslant n}\max_{z\in K\cap L^\perp}\|z\|\), Lemma~\ref{lem-section_width} gives the bound \(d_n(K)\leqslant(1+\sqrt n)d^n(K)\) (see~\cite[\S~2.10]{Pie87} for general relations between widths). The order of growth of the factor \(1+\sqrt n\) is optimal. Indeed, for \(n\geqslant1\) and \(K=B_{\ell_1^{2n}}\subset\ell_2^{2n}\), Theorem~\ref{th-l2} gives \(d_n(K)\geqslant1/\sqrt2\). By duality of widths~\cite{Pie87} and Kashin's bound~\cite[Theorem~1]{Kas77},
\[
d^n(K)=d_n(B_{\ell_2^{2n}},\ell_\infty^{2n})\leqslant C n^{-1/2}.
\]
\end{remark}

\begin{lemma}\label{lem-extract}
Let \(K\) be a compact convex set symmetric about the origin. If \(d_k(K)>0\) for all \(k\), then there exist an orthonormal system \(\{e_k\}_{k\geqslant0}\) and numbers \(r_k\geqslant d_k(K)/(1+\sqrt k)\) such that \(r_k e_k\in K\).
\end{lemma}

\begin{proof}
The vectors \(e_k\) and the numbers \(r_k\) are constructed inductively. Suppose that \(e_0, \dots, e_{m-1}\) have already been constructed, and let \(L_m = \operatorname{span}\{e_0, \dots, e_{m-1}\}\) (\(L_0 = \{0\}\)). By Lemma~\ref{lem-section_width}, there exists \(z_m \in K \cap L_m^\perp\) with \(\|z_m\| \geqslant d_m(K)/(1+\sqrt m) > 0\). Setting \(r_m = \|z_m\|\) and \(e_m = z_m/r_m\), we obtain \(r_m e_m \in K\) and \(e_m \perp L_m\).
\end{proof}

\begin{lemma}\label{lem-absconv_inner}
Let \(A\subset X\) be a nonempty compact set and put \(K=\overline{\operatorname{absconv}}(A)\). For every finite-dimensional subspace \(L\subset X\),
\[
\max_{x\in K}\|P_{L^\perp}x\|=\max_{a\in A}\|P_{L^\perp}a\|.
\]
For a unit vector \(e\in X\) and \(r>0\), condition~(iii) of Definition~\ref{def-spread} for the pair \((e,r)\) is equivalent to the inclusion \(re\in K\).
\end{lemma}

\begin{proof}
By Mazur's theorem, \(K\) is compact. The seminorm \(x\mapsto\|P_{L^\perp}x\|\) is continuous, convex and even, so its maximum on \(A\) bounds its values on \(K\). The inclusion \(A\subset K\) gives equality. If \(re\in K\), then condition~(iii) holds at a point where the maximum on \(A\) is attained.

Conversely, suppose that condition~(iii) holds. For a unit vector \(v\in X\) and \(\varepsilon>0\), let \(L\) be the linear span in \(v^\perp\) of a finite \(\varepsilon\)-net of the compact set \(P_{v^\perp}A\). Then
\[
r|\langle e,v\rangle|
\leqslant r\|P_{L^\perp}e\|
\leqslant\max_{a\in A}\|P_{L^\perp}a\|
\leqslant\max_{a\in A}|\langle a,v\rangle|+\varepsilon.
\]
Since \(\max_{a\in A}|\langle a,v\rangle|=\sup_{z\in K}\langle z,v\rangle\), letting \(\varepsilon\to0\) gives \(\langle re,v\rangle\leqslant\sup_{z\in K}\langle z,v\rangle\) for every unit \(v\). The separation theorem implies \(re\in K\).
\end{proof}

\begin{theorem}\label{th-main_dichotomy}
Let \(A\subset X\) be a compact set with \(0\in A\), and suppose that \(d_k(A)>0\) for all \(k\). Then there exist an orthonormal system \(\{e_k\}_{k\geqslant0}\) and numbers \(r_k\geqslant d_k(A)/(1+\sqrt k)\) for which condition~(iii) of Definition~\ref{def-spread} holds. In particular, if, for some \(C,\kappa>0\), we have \(d_k(A)\geqslant Ce^{-\kappa k}\) for all \(k\), then \(A\) is a strengthening set with parameters \((Cc_\varepsilon,\kappa+\varepsilon)\) for every \(\varepsilon>0\), where \(c_\varepsilon=\inf_{k\geqslant0}e^{\varepsilon k}/(1+\sqrt k)>0\).
\end{theorem}

\begin{proof}
By Lemma~\ref{lem-absconv_inner}, the set \(K=\overline{\operatorname{absconv}}(A)\) is compact, convex and symmetric about the origin, and \(d_k(K)=d_k(A)\). Lemma~\ref{lem-extract} gives an orthonormal system and points \(r_ke_k\in K\) with \(r_k\geqslant d_k(A)/(1+\sqrt k)\), and Lemma~\ref{lem-absconv_inner} gives condition~(iii). If \(d_k(A)\geqslant Ce^{-\kappa k}\), this implies \(r_k\geqslant Cc_\varepsilon e^{-(\kappa+\varepsilon)k}\). The number \(c_\varepsilon\) is positive, since \(e^{\varepsilon k}/(1+\sqrt k)\to\infty\).
\end{proof}

\begin{corollary}\label{cor-dichotomy_equiv}
Let \(A \subset X\) be a compact set with \(0 \in A\). The following conditions are equivalent:
\begin{enumerate}
\item[(a)] \(A\) is a strengthening set with some parameters \((C, \kappa)\);
\item[(b)] there exist \(C', \kappa' > 0\) such that \(d_n(A) \geqslant C' e^{-\kappa' n}\) for all \(n \in \mathbb{Z}_+\).
\end{enumerate}
The following dichotomy holds: either \(A\) is a strengthening set with some parameters, or
\[
\limsup_{n \to \infty} \frac{-\ln d_n(A)}{n} = +\infty \qquad (-\ln 0 := +\infty).
\]
\end{corollary}

\begin{proof}
Condition~(a) with parameters \((C, \kappa)\) implies~(b) with \(C' = Ca(\kappa)\) and \(\kappa' = \kappa\) (Theorem~\ref{th-rootfree}), and~(b) implies~(a) with parameters \((C' c_\varepsilon, \kappa'+\varepsilon)\) for every \(\varepsilon > 0\) (Theorem~\ref{th-main_dichotomy}). It remains to check that~(b) is equivalent to finiteness of the limit superior. Condition~(b) gives \(-\ln d_n(A)/n \leqslant \kappa' - (\ln C')/n\), so the limit superior is at most \(\kappa'\). Conversely, suppose that the limit superior is \(S < \infty\) and \(\kappa' > \max\{S, 0\}\). Then \(d_n(A) \geqslant e^{-\kappa' n}\) for \(n \geqslant N\) with some \(N\geqslant1\). Since \(d_N(A)>0\), we have \(d_n(A) \geqslant d_N(A) > 0\) for \(n < N\) by monotonicity of the widths in \(n\), and~(b) holds with \(C' = \min\{1, \min_{n<N} d_n(A) e^{\kappa' n}\}\).
\end{proof}

\section{Applications to dissipative systems}\label{sec-appl}

\subsection{A diagonal dissipative system}\label{sec-pde}

Let \(\{\varphi_k\}_{k\geqslant0}\) be an orthonormal basis of \(H=L^2(\mathbb T^d)\). Consider the diagonal system
\begin{equation}\label{eq-rd_system}
\dot u_k=u_k-c_k^{-2}u_k^3, \qquad k\geqslant0,
\end{equation}
where \(u_k=\langle u,\varphi_k\rangle\) and \(c_k=Ce^{-\gamma k}\), \(C,\gamma>0\). The system is understood coordinatewise, and its solutions are regarded as continuous \(H\)-valued functions.

The solutions of the scalar equations exist for all \(t\geqslant0\), since \(|u_k(t)|\leqslant e^t|u_k(0)|\). The right-hand side \(f_k(\xi)=\xi-c_k^{-2}\xi^3\) satisfies
\[
(f_k(\xi)-f_k(\zeta))(\xi-\zeta)\leqslant |\xi-\zeta|^2.
\]
Thus the maps \(u_k(0)\mapsto u_k(t)\) are Lipschitz with constant \(e^t\) common to all \(k\). The solutions of the scalar equations define a continuous semiflow \(S(t)\) on \(H\). Continuity in time follows from the bound \(|u_k(t)|\leqslant e^t|u_k(0)|\) and the dominated convergence theorem, and joint continuity in time and initial data follows from the Lipschitz bound, which is uniform on bounded time intervals.

The parallelepiped \(\mathcal A=\{u\in H:|u_k|\leqslant c_k,\ k\geqslant0\}\) is compact, since \(\sum_k c_k^2<\infty\). For every \(t\geqslant0\) the map \(u_k(0)\mapsto u_k(t)\) takes \([-c_k,c_k]\) onto itself, so \(S(t)\mathcal A=\mathcal A\). For \(s=|u_k|>c_k\) we have
\[
\dot s=-c_k^{-2}s(s+c_k)(s-c_k)\leqslant-2(s-c_k).
\]
The distance to \(\mathcal A\) is computed coordinatewise, so
\[
\inf_{v\in\mathcal A}\|S(t)u-v\|
\leqslant e^{-2t}\inf_{v\in\mathcal A}\|u-v\|.
\]
Thus \(\mathcal A\) is the global attractor of system~\eqref{eq-rd_system}.

The set \(\mathcal A\) contains \(A_{\mathrm{par}}\) with \(e_k=\varphi_k\) and \(r_k=c_k\). Proposition~\ref{prop-box_exact} gives the lower bound, and the subspace \(\operatorname{span}\{\varphi_0,\ldots,\varphi_{n-1}\}\) gives the matching upper bound. Hence
\begin{equation}\label{eq-rd_result}
d_n(\mathcal A)
=\Bigl(\sum_{k\geqslant n}c_k^2\Bigr)^{1/2}
=\frac{Ce^{-\gamma n}}{\sqrt{1-e^{-2\gamma}}}.
\end{equation}

\subsection{A lower bound from stationary solutions}\label{sec-stationary}

\begin{proposition}\label{prop-stationary}
\begin{sloppypar}
Let \(\{\varphi_k\}_{k\geqslant0}\) be an orthonormal system in the phase space. Suppose that the global attractor \(\mathcal A\) contains the zero solution and a sequence of stationary solutions \(\{u_k\}_{k\geqslant0}\) of the form
\begin{equation}\label{eq-decomp}
u_k = c_k \varphi_k + w_k, \qquad w_k \perp \varphi_k,
\end{equation}
where the coefficients \(c_k\) satisfy the exponential lower bound
\begin{equation}\label{eq-amp2}
|c_k| \geqslant C_1 e^{-\gamma k} \qquad (C_1,\gamma>0),
\end{equation}
and the relative remainders are square summable:
\begin{equation}\label{eq-eta}
\sum_{k=0}^{\infty}\eta_k<\infty, \qquad \eta_k := \frac{\|w_k\|^2}{|c_k|^2}.
\end{equation}
Then there exists a constant \(C_*>0\) such that
\begin{equation}\label{eq-pde_result}
d_n(\mathcal A) \geqslant C_* e^{-\gamma n} \qquad \text{for all } n\in\mathbb Z_+.
\end{equation}
\end{sloppypar}
\end{proposition}

\begin{proof}
Choose \(m\) so that \(\sum_{k\geqslant m}\eta_k\leqslant1/4\). As in the proof of Theorem~\ref{th-almost_summable}, the operator
\[
Tx=x+\sum_{k\geqslant m}\langle x,\varphi_k\rangle\,\frac{w_k}{c_k}
\]
is invertible and \(\|T^{-1}\|\leqslant2\). Since \(T(c_k\varphi_k)=u_k\) for \(k\geqslant m\), the set \(B=\{0\}\cup\{c_k\varphi_k:k\geqslant m\}\) satisfies \(T(B)\subset\mathcal A\). The set \(B\) is bounded because \(|c_k|\leqslant\|u_k\|\), and by~\eqref{eq-amp2} it has orthogonal spread with parameters \((C_1e^{-\gamma m},\gamma)\). The bound for \(T^{-1}\) and Corollary~\ref{cor-orth_bound} give
\[
d_n(\mathcal A)\geqslant\tfrac12d_n(B)\geqslant\tfrac12 C_1e^{-\gamma m}a(\gamma)e^{-\gamma n}.
\]
\end{proof}

Condition~\eqref{eq-eta} holds, for example, if \(|c_k|\leqslant C_2e^{-\gamma k}\) and \(\|w_k\|\leqslant B|c_k|^2\) for all \(k\) with constants \(C_2,B>0\). Indeed, in this case \(\eta_k\leqslant B^2C_2^2e^{-2\gamma k}\).

\begin{remark}\label{rem-ezk}
In the example of Eden, Zelik and Kalantarov~\cite[Lemma~5.1, formulas~(5.14), (5.16), Theorem~5.4]{EZK13}, the orthogonal projection \(Q\mathcal A\) of the attractor onto \(H\) contains the points \(b_j e_j\), \(j\geqslant1\), where \(\{e_j\}\) is an orthonormal system in \(H\) and \(b_j\asymp j^{-3}(\ln(j+1))^{-3}\). Applying~\eqref{eq-l2M} to the first \(2n+1\) directions and noting that orthogonal projection does not increase widths, we obtain, in the norm of \(\mathbb R^2\oplus H\),
\[
d_n(\mathcal A)\geqslant d_n(Q\mathcal A)
\gtrsim (n+1)^{-3}(\ln(n+2))^{-3}.
\]
\end{remark}

\begin{remark}\label{rem-allen_cahn}
Without giving a detailed proof here, we note that the method applies to the global attractor \(\mathcal A_\nu\) of the Allen--Cahn equation on the interval \([0,2\pi]\) with periodic boundary conditions,
\[
u_t=\nu u_{xx}+u-u^3,\qquad 0<\nu<1.
\]
In \(L^2(0,2\pi)\), the method gives the bound
\[
d_n(\mathcal A_\nu)\geqslant C_\nu e^{-\kappa_\nu n},
\qquad n\in\mathbb Z_+,
\]
where the constants \(C_\nu,\kappa_\nu>0\) depend only on \(\nu\). It suffices to use the spatial translates of a single nontrivial stationary solution. Its expansion in odd Fourier harmonics has coefficients satisfying an exponential lower bound. The corresponding formulas are given in~\cite[formulas~(22.11.1), (22.13.13)]{DLMF}. Isolating individual harmonics gives line segments in the closed absolutely convex hull of the attractor, parallel to the vectors of an orthonormal system, with half-lengths satisfying an exponential lower bound. By Lemma~\ref{lem-absconv_inner}, this gives the strengthening property, and Theorem~\ref{th-l2} then yields the bound.
\end{remark}

\textbf{AI disclosure statement.} ChatGPT (OpenAI) was used as a discussion partner to explore further implications of the strengthening method. It also assisted in applying the method to concrete partial differential equations, including the periodic Allen–Cahn equation. The author takes full responsibility for the mathematical content and the final text.


\begin{thebibliography}{99}

\bibitem{Kol36}
A.~Kolmogoroff,
\emph{\"Uber die beste Ann\"aherung von Funktionen einer gegebenen Funktionenklasse},
Ann. of Math. (2) \textbf{37}:1 (1936), 107--110.

\bibitem{Tik60}
V.\,M.~Tikhomirov,
\emph{Diameters of sets in function spaces and the theory of best approximations},
Russian Math. Surveys \textbf{15}:3 (1960), 75--111.

\bibitem{Pin85}
A.~Pinkus,
\emph{$n$-Widths in Approximation Theory},
Springer-Verlag, Berlin, 1985.

\bibitem{KolTik59}
A.\,N.~Kolmogorov and V.\,M.~Tikhomirov,
\emph{$\varepsilon$-entropy and $\varepsilon$-capacity of sets in function spaces},
Uspekhi Mat. Nauk \textbf{14}:2 (1959), 3--86.

\bibitem{BabVis83}
A.\,V.~Babin and M.\,I.~Vishik,
\emph{Attractors of partial differential evolution equations and estimates of their dimension},
Russian Math. Surveys \textbf{38}:4 (1983), 151--213.

\bibitem{Tem97}
R.~Temam,
\emph{Infinite-Dimensional Dynamical Systems in Mechanics and Physics},
Springer-Verlag, New York, 1997.

\bibitem{CheVis02}
V.\,V.~Chepyzhov and M.\,I.~Vishik,
\emph{Attractors for Equations of Mathematical Physics},
Amer. Math. Soc., Providence, RI, 2002.

\bibitem{Zel23}
S.\,V.~Zelik,
\emph{Attractors. Then and now},
Russian Math. Surveys \textbf{78}:4 (2023), 635--777.

\bibitem{EZK13}
A.~Eden, S.\,V.~Zelik and V.\,K.~Kalantarov,
\emph{Counterexamples to regularity of Ma\~n\'e projections in the theory of attractors},
Russian Math. Surveys \textbf{68}:2 (2013), 199--226.

\bibitem{Coh14}
A.~Cohen,
\emph{Estimating the $n$-width of solution manifolds of parametric PDE's},
slides of a talk, Chicheley Hall, June~30, 2014.

\bibitem{Ism68}
R.\,S.~Ismagilov,
\emph{On $n$-dimensional diameters of compacts in a Hilbert space},
Funct. Anal. Appl. \textbf{2}:2 (1968), 125--132.

\bibitem{Sof69}
L.\,B.~Sofman,
\emph{Diameters of octahedra},
Math. Notes \textbf{5}:4 (1969), 258--262.

\bibitem{Henk12}
M.~Henk,
\emph{L\"owner--John ellipsoids},
In Optimization Stories. Doc. Math. Extra Vol. ISMP (2012), 95--106.

\bibitem{Pie87}
A.~Pietsch,
\emph{Eigenvalues and $s$-Numbers},
Cambridge University Press, Cambridge, 1987.

\bibitem{Kas77}
B.\,S.~Kashin,
\emph{Diameters of some finite-dimensional sets and classes of smooth functions},
Math. USSR-Izv. \textbf{11}:2 (1977), 317--333.

\bibitem{DLMF}
F.\,W.\,J.~Olver, A.\,B.~Olde Daalhuis, D.\,W.~Lozier,
B.\,I.~Schneider, R.\,F.~Boisvert, C.\,W.~Clark, B.\,R.~Miller,
B.\,V.~Saunders, H.\,S.~Cohl and M.\,A.~McClain (eds.),
\emph{NIST Digital Library of Mathematical Functions},
Release~1.2.8 of September~15, 2026.

\end{thebibliography}
\end{document}